\documentclass[final,leqno,onefignum,onetabnum,dvipsnames]{siamltex1213}
\usepackage[normalem]{ulem}
\usepackage{framed}
\usepackage{amssymb}
\usepackage{graphicx}
\usepackage[outdir=./]{epstopdf}
\usepackage{float}
\usepackage{caption}
\usepackage{sidecap}
\usepackage{wrapfig}
\usepackage{color}
\usepackage{mathptmx}
\RequirePackage{ifthen}
\usepackage{times}
\usepackage{amsmath}
\usepackage{algpseudocode}
\usepackage{enumerate}
\usepackage{color}
\usepackage{xcolor}
\usepackage{mathtools}
\DeclarePairedDelimiter\ceil{\lceil}{\rceil}
\DeclarePairedDelimiter\abs{\lvert}{\rvert}
\DeclarePairedDelimiter\norm{\lVert}{\rVert}

\DeclareMathOperator{\Cov}{Cov}
\DeclareMathOperator{\Tr}{Tr}

\title{A Seamless Multilevel Ensemble Transform Particle Filter \footnote{Preprint submitted to arXiv.}} 

\author{A. Gregory \footnotemark[2] and C. J. Cotter \footnotemark[2]} 

\begin{document}
\maketitle

\renewcommand{\thefootnote}{\fnsymbol{footnote}}

\footnotetext[2]{Department of Mathematics, Imperial College London, Exhibition Road, London SW7 2AZ, UK (a.gregory14@imperial.ac.uk). Alastair Gregory was supported by the Science and Solutions to a Changing Planet DTP and the Natural Environmental Research Council.}

\renewcommand{\thefootnote}{\arabic{footnote}}

\slugger{mms}{xxxx}{xx}{x}{x--x}%slugger should be set to mms, siap, sicomp, sicon, sidma, sima, simax, sinum, siopt, sisc, or sirev

\begin{abstract}
This paper presents a seamless algorithm for the application of the
  multilevel Monte Carlo (MLMC) method to the ensemble transform
  particle filter (ETPF). The algorithm uses a combination of optimal
  coupling transformations between coarse and fine ensembles in
  difference estimators within a multilevel framework, to minimise
  estimator variance. It differs from that of Gregory et al. (2016) in
  that strong coupling between the coarse and fine ensembles is
  seamlessly maintained during all stages of the assimilation algorithm, instead of
  using independent transformations to equal weights followed by
  recoupling with an assignment problem.  This modification is found
  to lead to an increased rate in variance decay between coarse and
  fine ensembles with level in the hierarchy, a key component of
  MLMC. This offers the potential for greater computational cost
  reductions. This is shown, alongside evidence of asymptotic
  consistency, in numerical examples.
\end{abstract}



\pagestyle{myheadings}
\thispagestyle{plain}
\markboth{A. GREGORY AND C. J. COTTER}{A SEAMLESS MULTILEVEL ENSEMBLE PARTICLE TRANSFORM FILTER}

\section{Introduction}

Particle filters and some parametric filters use ensembles to represent a posterior measure given a partially observed random dynamical system. Recently, a number of studies \cite{Gregory, Hoel, Jasra} have considered the extension of the multilevel Monte Carlo method to filtering problems \cite{Giles}.  This involves using a telescoping sum of Monte Carlo estimators of solutions to systems discretized with varying levels of resolution. The method has also been applied to other Bayesian inference problems, such as static parameter estimation \cite{JasraStatic}.

The applications of MLMC  have varied greatly, with the only major difference being  how to `couple' random samples of the system from two consecutive levels to achieve variance reduction in the hierarchy of Monte Carlo estimators. For nonlinear filtering, \cite{Jasra} uses a coupled resampling algorithm for a particle filter (multilevel particle filter) whereas \cite{Gregory} used an optimal transportation coupling algorithm for an ensemble transform particle filter (ETPF) \cite{Reich} leading to the multilevel ETPF (MLETPF). Both of these algorithms returned computational cost reductions for a fixed order of accuracy from that of their standard, single level counterparts.

Since then, the optimal transport algorithm has been investigated in further detail in \cite{Sen} with the aim of the coupling of multiple particle filters. This study also considered the use of the iterative Sinkhorn algorithm, as well as exploiting sparse optimal transport matrices, to solve these optimal transportation problems, resulting in the reduction of their computational cost.
The iterative Sinkhorn algorithm has also been implemented to reduce the computational cost of the optimal transportation problem in a second-order accurate ETPF \cite{DeWiljes}.

This paper will concentrate on improving the algorithm in \cite{Gregory} to avoid using an assignment problem \cite{Hungarian}, instead using a `seamless' combination of three optimal transportation problems. 
The benefit of applying this modification is that, for certain numerical cases considered in this paper, the rate of variance decay between samples with consecutive levels of increasing resolution is improved from that in \cite{Gregory}; this leads to a greater rate in the reduction of computational cost for the propagation of ensembles in the MLETPF. 
In the notation of \cite{Giles}, we observe a variance decay scaling close to $\beta \approx 2$ (quadratic decay) in multivariate problems rather than one that corresponds to $\beta\approx 1$ (linear decay) which was observed previously. In other cases, including ones where the problem is univariate or localisation is used, both the proposed and previous algorithms scale at the optimal quadratic rate. However the proposed algorithm is still found to return estimators with a smaller magnitude of variance and thus are less costly for a fixed accuracy.

\iffalse

Second, the modification frames the problem as three optimal transportation problems; this could allow the algorithm to benefit from recent developments in iterative optimal transportation research \cite{Ravi, Li} which has the potential to reduce the computational cost of the optimal transport resampling step in the ETPF. So far, this has only been trialed in a second-order version of the ETPF \cite{DeWiljes}. This is a point that has been left for consideration in a later study and is not investigated in this paper. It must be noted that this will likely be a non-trivial extension for the multilevel framework and poses a challenge for such a study.

\fi

The paper proceeds as follows. In Section 2, an introduction to Bayesian filters will be given, and then a new `seamless' agorithm for the optimal transport resampling step in a variant of the particle filter, the ETPF, will be explained. Numerical examples will provide a proof of concept in Section 3. Finally a summary and discussion concludes the paper in Section 4.

\section{Multilevel Filtering Algorithms}

\subsection{Particle Filters for Bayesian Data Assimilation Problems}

\bigskip

\label{subsec:subsec2_1}

We consider a Bayesian data assimilation problem for a state vector $X\in
\mathbb{R}^d$, with discrete time model equation
\begin{equation}
  X^{n+1} = M(X^n;\omega),
\end{equation}
where $\omega$ denotes any stochastic contributions in $M$, and where
$M$ may involve the application of several timesteps with the index
$n$ denoting the observation times. We assume that at each observation
time, $t_{n}$ (where the interval between two observation times, $\Delta t = t_{n+1} - t_{n}$, is assumed fixed), we have noisy observations $Y^n\in \mathbb{R}^m$ given by
\begin{equation}
  Y^n = h(X^n) + \eta^n,
\end{equation}
where $h:\mathbb{R}^d\to\mathbb{R}^m$ is a (nonlinear) observation
operator, and $\eta^n\in \mathbb{R}^m$ is a random variable
representing measurement error that is assumed in this paper to be drawn from $N(0, R)$, $R \in \mathbb{R}^{m \times m}$, so that $Y^n$ has probability density
\begin{equation}
  \pi(Y^n;X^n)  = \frac{1}{\sqrt{2 \pi R}} \exp\left(-\frac{1}{2}\left(Y^{n}-h(X^{n})\right)^TR^{-1}\left(Y^{n}-h(X^{n})\right)\right),
\end{equation}
denoted the likelihood function. A  filter is an estimator for
the conditional probability density
\begin{equation}
  \pi(X^n|Y^0,\ldots,Y^n) 
\end{equation}
for the state variable $X^n$, given the observations until time index
$n$ and assuming a given initial probability density $\pi(X_0)$. The classical particle
 filter uses an ensemble $\{X^{i, n}\}_{i=1}^N$ of model states (here $i$ denotes the realisation of model state and $n$ denotes the time index)
together with a set of scalar weights $\{w^{i, n}\}_{i=1}^N$, so that
expectations over the conditional probability density may be consistently estimated
using the formula
\begin{equation}
  E[X^n;Y^0,\ldots,Y^n] \approx
  \hat{E}^{N}[X^n;Y^0,\ldots,Y^n] = \sum_{i=1}^N w^{i, n}X^{i, n}.
\end{equation}
A particle filter can be constructed iteratively using three stages.
\begin{enumerate}
\item {\bfseries Model propagation}:
  \begin{equation}
    X^{i, n} = M(X^{i, n-1};\omega^i), i=1,\ldots,N.
  \end{equation}
        (Here $\omega^i$ represents an independent realisation
    of the random variable $\omega$ for each ensemble member $X^{i, n}$.)
  \item {\bfseries Bayesian weight update}:
    \begin{equation}
      w^{i, n} = \frac{1}{Z}w^{i, n-1}\pi(Y^n|X^{i, n}),
    \end{equation}
    where $Z$ is a normalisation constant chosen so that $\sum_{i=1}^Nw^{i, n}=1$.
  \item {\bfseries Resampling / transformation}: A (possibly stochastic) ensemble
    transformation $\{X^{i, n},w^{i, n}\}_{i=1}^N\mapsto
    \{\tilde{X}^{i, n},1/N\}_{i=1}^N$. There are deterministic cases of this transformation \cite{Cotter} that preserve the filter estimate
    \begin{equation}
    \hat{E}^{N}[X^n;Y^0,\ldots,Y^n] = \sum_{i=1}^N w^{i, n}X^{i, n} = \frac{1}{N} \sum^{N}_{i=1}\tilde{X}^{i, n} ,
    \end{equation}
    or in the stochastic transformation case (random resampling) provide an unbiased estimate
    \begin{equation}
    E\Bigg[ \frac{1}{N} \sum^{N}_{i=1}\tilde{X}^{i, n}\Bigg] = \sum_{i=1}^N w^{i, n}X^{i, n} .
    \end{equation}
\end{enumerate}
Specific particle filters are defined through the choice of
resampling/transformation method, which is necessary to avoid degenerate
weights \cite{Cappe}. In particular, \cite{Doucet} provides a background on a range of resampling techniques, including adaptive resampling, and their effects on the particle filter. Typically, random resampling methods add variance to the filtering estimator. There is literature exploring the rigorous theory behind the impact of this on the particle filter, such as in \cite{Chopin}. Resampling ideas like residual resampling try to reduce this additional variance. The deterministic resampling/transformation described below, which was proposed in \cite{Reich}, is designed to minimise this additional variance.

\subsection{Optimal Coupling Between Discrete Probability Distributions}
The key tool in the type of resampling considered in this paper is an optimal coupling between discrete
probability distributions. First, consider a pair $(Z_1, Z_2)$ of discrete
random variables where $Z_{1}$ take values in the set $\{Z_1^i\}_{i=1}^N$, with $Z_1^i \in \mathbb{R}^{d}$,
and $Z_2$ takes values in the set $\{Z_2^i\}_{i=1}^N$, with $Z_2^i \in \mathbb{R}^{d}$. We wish to obtain a
joint probability distribution
\begin{equation}
  P(Z_1=Z_1^i,Z_2=Z_2^j) = D_{i,j}, \quad i,j=1,\ldots,N,
\end{equation}
with marginal probabilities
\begin{equation}
\begin{split}
  &P(Z_1=Z_1^i) = \sum_{j=1}^N D_{i,j} = p^i, \, i=1,\ldots,N, \\
  \quad & P(Z_2=Z_2^j) = \sum_{i=1}^N D_{i,j} = q^j, \, j=1,\ldots,N.
\end{split}
\end{equation}
For one to minimise the variance of the difference between $Z_{1}$ and $Z_{2}$, one requires the joint probability matrix $D_{i,j}$ that maximises the
covariance $\Cov_{D}(Z_1,Z_2)$ subject to these marginals. We note
\begin{equation}
E_{D}[\norm*{Z_{1}-Z_{2}}^{2}]=E[\norm*{Z_{1}}^2]+E[\norm*{Z_{2}}^2]-2E[Z_{1}]^{T}E[Z_{2}]-2\Tr\left(\Cov_{D}[Z_{1},Z_{2}]\right)
\end{equation}
where $E_{D}$  and $\Cov_{D}$ are the expectation and covariance with respect to the coupling probability matrix $D$. Thus maximising the covariance $\Cov_{D}(Z_1,Z_2)$ is
equivalent to minimising the cost function (the discrete Wasserstein distance)
\begin{equation}
  J(D) = \sum_{i=1}^N\sum_{j=1}^N D_{i,j}\norm*{Z_1^i - Z_2^j}^2,
\end{equation}
over $D_{i,j}$, subject to the constraints
\begin{equation}
  \sum_{i=1}^N D_{i,j} = q^j,
  \quad \sum_{j=1}^N D_{i,j} = p^i.
\end{equation}
The solution to this optimisation problem can be found by linear
programming. The computational cost of this can be $\mathcal{O}(N\log N)$ if the dimension, $d$, is one, as it can be solved by a cheap algorithm \cite{Cotter}, but is $\mathcal{O}(N^{3}\log N)$ in multiple dimensions. For cheap forward models, this can be a computational bottleneck. Localisation, which is described in the next section, can provide some relief in this issue. As described in \cite{Cotter}, the resulting matrix
$D_{i,j}$ is typically very sparse, and in the one-dimensional case, has a maximum of two non-zeros per row. This optimal coupling can be used to resample weighted
ensembles.

Let the ensembles $\{Z_1^i\}_{i=1}^N$ and $\{Z_2^i\}_{i=1}^N$ be
drawn from random variables $Z_{1}$ and $Z_{2}$ with measures $\mu_{Z_{1}}$ and $\mu_{Z_{2}}$ respectively. Furthermore suppose they are assigned weights $\{p^i\}_{i=1}^N$ and $\{q^i\}_{i=1}^N$ respectively by importance sampling \emph{via} a third probability measure
$\mu_{Z_3}$, all absolutely continuous with respect to the Lesbegue
measure. We wish to resample $\{Z_2^i,q^i\}_{i=1}^N$ to 
obtain a new ensemble $\{\tilde{Z}_2^i,p_i\}_{i=1}^N$, such that
\begin{equation}
  E\left[\sum_{i=1}^N p^i\tilde{Z}_2^i\right] = \sum_{i=1}^N q^i Z_2^i,
\end{equation}
where the expectation is taken over the random variables used in the
resampling. This can be achieved by choosing $\tilde{Z}_2^i = Z_2^j$ with
probability $D_{i,j}/p^i$. Then
\begin{equation}
  E\left[\sum_{i=1}^N p^i\tilde{Z}_2^i\right] = \sum_{i=1}^N p^i\sum_{j=1}^N \left(\frac{D_{i,j}}{p^i}\right)Z_2^j
  = \sum_{j=1}^N \left(\sum_{i=1}^N D_{i,j}\right)Z_2^j = \sum_{j=1}^N q^j Z_2^j,
\end{equation}
as required.\\
As an alternative to this random resampling, one can also take the
deterministic transformation
\begin{equation}
  \tilde{Z}_2^i = \sum_{j=1}^N\left(\frac{D_{i,j}}{p^i}\right)Z_2^j.
\label{equation:transform}
\end{equation}
In \cite{Cotter}, the following was proved. If we define the set $\{\tilde{Z}_2^i\}_{i=1}^{N}$, computed via the transform in (\ref{equation:transform}), then the maps $\Psi_N:\{Z_1^i\}_{i=1}^N\to
\{\tilde{Z}_2^i\}_{i=1}^N$, given by
$$\tilde{Z}_{2}^{i} = \Psi_N(Z_{1}^{i}), \qquad i=1,...,N, $$ 
weakly converge to a map $\Psi:\mathbb{R}^d\to\mathbb{R}^d$ as $N\to \infty$. In addition to this, the random variable given by $Z_{2} = \Psi(Z_{1})$ has probability measure $\mu_{Z_{2}}$.
Hence, for arbitrary test functionals $g$, we can obtain weak convergence for estimates to $E[g(Z_{2})]$, and
\begin{equation}
\sum_{i=1}^N p^i g(\tilde{Z}_2^i) \to \sum_{i=1}^N q^ig(Z_2^i),
\end{equation}
as $N\to \infty$.

The optimal coupling can be used to define a deterministic ensemble
transform which was introduced in \cite{Reich} to replace
a weighted ensemble $\{Z_1^i,w^i\}_{i=1}^N$ by an equally weighted
ensemble $\{\tilde{Z}_1^i,1/N\}_{i=1}^N$. We take $Z_2^i:=Z_1^i$, $i=1,\ldots,N$,
in the above coupling scheme and compute
\begin{equation}
  \tilde{Z}_1^j = N \sum_{i=1} D_{i,j}Z_1^i.
\end{equation}
This gives a consistent ensemble with convergent statistics as $N\to
\infty$. This is the basis of the Ensemble Transform Particle Filter
(ETPF) of \cite{Reich}, where this transform is used as a
resampling step. It would be desirable to establish $L_{2}$ rates of convergence for this transformation; this could allow one to establish rigorous analysis of the effectiveness of the multilevel Monte Carlo framework which will be described later.

\subsection{Localisation}

Localisation is useful in cases where the dimension $d$ is large or one is dealing with spatially extended systems. A multilevel version of the Ensemble Kalman Filter (EnKF) for these types of system has recently been proposed \cite{Chernov}. Localisation in the ETPF can be used to simplify the optimal transportation problems used above, and in some cases, reduce the computational cost of them from $\mathcal{O}(N^3logN)$ to $\mathcal{O}(dNlogN)$. Using a ``localisation radius", $r_{loc}$, one can solve a separate optimal transportation problem on each individual component of the random variables $Z_{1}$ and $Z_{2}$, with the cost function for the $k$'th component of $Z_{1}$ and $Z_{2}$ (denoted $Z_{1}(k)$ and $Z_{2}(k)$ respectively) being
\begin{equation}
 J(D) = \sum^{d}_{m=1}\sum_{i=1}^N\sum_{j=1}^N C_{k,m}D_{i,j}\left(Z_1^i(k) - Z_2^j(k)\right)^2,
\end{equation}
where $C$ is some localisation matrix. This could be of the form
\begin{equation}
C_{k,m}=
\begin{cases}
1-\frac{1}{2}\left(\frac{s}{r_{loc}}\right), & \quad s \leq 2,\\
0, & \quad \text{otherwise},
\end{cases} 
\end{equation}
where $s=\min\big\{\abs{k-m-N},\abs{k-m},\abs{k-m+N}\big\}$. In the case when $r_{loc}=0$,
each optimal transportation problem is one-dimensional (can be solved in $\mathcal{O}(NlogN)$), and the problems for all individual components can be solved in $\mathcal{O}(dNlogN)$. This case will be referred to as the fully localised case. It must be noted however, that localisation adds an additional bias into the ETPF estimate.
For more information on this localisation scheme that can be used alongside the ETPF refer to \cite{Cheng}.

Localisation in the ETPF is important in cases where the optimal transport algorithmic cost is dominative over the model cost of the problem. In these cases, the multilevel framework (which can return computational cost reductions in the cost of the forward model / propagation of ensembles) would prove to have a negligible impact on the overall cost of the filter. Thus localisation should be used in those cases. Localisation in the above manner is used in one of the numerical examples at the end of this paper.
\iffalse
Therefore when considering the effectiveness of the proposed multilevel version of the ETPF, it is assumed that the forward model cost dominates over that of solving the optimal transport problem, either with full localisation or without (e.g. in some spatially extended systems).
\fi

\subsection{A Multilevel Filter}

\label{subsec:ml_filter}

The recent emergence of the multilevel Monte Carlo method poses the question of how to apply the telescoping sum framework to filters \cite{Gregory, Jasra, Hoel}.  A multilevel filter estimate makes use of a hierarchy of models $M^l$, $l=0,\ldots,L$, where $l=0$ denotes the coarsest,
cheapest and least accurate model and $l=L$ denotes the finest, most
expensive and most accurate model (in this paper we only consider models
with different time-step sizes so that all state vectors have the same
dimension $d$). We assume that each model has a time-step size $h_l$
and that the model cost scales as $\mathcal{O}(h_l^{-\gamma})$ for
$\gamma>0$. A multilevel filter uses a coarse level ensemble
$\{X^{i, n}_0,w^{i, n}_0\}_{i=1}^{N_0}$ (of i.i.d. model realisations $X^{i, n}_0$ where the subscript denotes the level of resolution of the realisation and $n$ is the time-index) together with a hierarchy of
pairs of ensembles
$\{\hat{X}^{i, n}_{l-1},\hat{w}^{i, n}_{l-1},X^{i, n}_{l},w^{i, n}_{l}\}_{i=1}^{N_{l}}$,
for $l=1,\ldots,L$.
Here $\hat{X}_l^{i, n} \sim X_l^n$ but the samples $\hat{X}_{l}^{i, n}$, with importance weights $\hat{w}_{l}^{i, n}$, are independent to the samples $X_{l}^{i, n}$. Thus expectations may be estimated using the telescoping sum
formula
\begin{equation}
  E[X_{L}^n] \approx \hat{E}^{N_{0},...,N_{L}}[X_{L}^{n}], \qquad \hat{E}^{N_{0},...,N_{L}}[X_{L}^{n}]:\hat{E}^{N_{0}}[X_{0}^n]
  + \sum_{l=1}^L \hat{E}^{N_{l}}[X_{l}^n-X_{l-1}^n],
\end{equation}
where we define the coarse estimator
\begin{equation}
  \hat{E}^{N_{0}}[X_{0}^n] = \sum_{i=1}^{N_0} w_{0}^{i, n} X_{0}^{i, n},
\end{equation}
and the difference estimators
\begin{equation}
  \hat{E}^{N_{l}}[X_{l}^n-X_{l-1}^n] = \sum_{i=1}^{N_l} \left(w_{l}^{i, n} X_{l}^{i, n}
  - \hat{w}_{l-1}^{i, n} \hat{X}_{l-1}^{i, n}\right).
\end{equation}
If the two independent estimators $\hat{E}^{N_{l}}[X_{l-1}^n]$ and $\hat{E}^{N_{l-1}}[X_{l-1}^n]$ are both asymptotically consistent estimators of $E[X_{l-1}^n]$ then recursively the telescoping sum estimator is consistent with $E[X_{L}^n]$. In filtering one is also interested in estimating higher moments of $X_{L}^{n}$ and thus being able to asymptotically consistently estimate $E[g(X_{L}^{n})]$ using this telescoping sum is also important. Despite the weak convergence of these estimators mentioned throughout this paper, establishing $L_{2}$ convergence rates for them remains an open problem.

The three iterative filtering steps defined in Section \ref{subsec:subsec2_1} are applied to this
set of $2L + 1$ ensembles, noting that: (a) the random sampling of the
initial conditions, as well as the stochastic terms in $M$, must be
drawn independently for the coarse estimator and each difference
estimator; and (b) the correlation between the coarse and fine level
ensembles in each difference estimator must be as high as possible. In
general for multilevel Monte Carlo methods, the latter is achieved by
using the same initial conditions for the coarse and fine ensembles in
each difference estimator, as well as the same realisations of the
stochastic terms in $M$. 

If the correlation between the coarse and fine level ensembles in each difference estimator
is sufficiently great then we satisfy a key condition of the multilevel Monte Carlo framework. This is that the variance of the difference between pairs of samples in the weighted coarse and fine ensembles, $\mathbb{V}_{l}$, decreases asymptotically with increasing $l$. When this condition is achieved, and we allow the sample sizes of each of the independent estimators, $N_{l}$, to decrease asymptotically with $l$ at a certain rate, then generally speaking, model computational cost reductions from a standard single level filter with fixed accuracy can be achieved. This is because each independent estimator balances estimator variance (determined by $V_{l}/N_{l}$) and discretization bias; expensive estimators with small discretization bias have less samples and cheap estimators with large discretization bias have more samples.

For the standard multilevel Monte Carlo methodology in \cite{Giles}, an algorithm was presented to find the optimal values of $L$ and $N_{l}$ that produce the greatest computational cost reductions from the Monte Carlo single level counterpart with the same accuracy. For more analysis and accompanying theory, turn to the review article \cite{GilesReview}. This article also gives a brief guide to the range of applications of multilevel Monte Carlo.

For multilevel filtering algorithms, the
challenge is to find a resampling/transformation strategy that keeps
the coarse and fine ensembles correlated in each difference estimator after each assimilation step.

\subsection{The Multilevel Ensemble Transform Particle Filter}

In \cite{Gregory}, we proposed an algorithm to address this issue. After
updating the weights of each pair of ensembles, we took the following steps to apply the multilevel Monte Carlo method to the ETPF, creating the multilevel ensemble transform particle filter (MLETPF).

\bigskip

\begin{enumerate}
\item {\bfseries Separately transform each ensemble}. 

  We independently apply an ensemble transform to both the coarse and fine ensemble in each difference estimator, following the approach
  of \cite{Reich}. This transform was described earlier in the paper. For both the coarse and fine weighted ensembles in each difference estimator, with $N_{F}$ particles in each, we seek a joint probability (coupling) matrix $D_{i,j}$ between $\{X^{i,n},w^{i,n}\}_{i=1}^{N_{F}}$ and the evenly weighted ensemble $\{\tilde{X}^{i,n}\}_{i=1}^{N_{F}}$. In particular we desire the $D_{i,j}$, $i=1,...,N_{F}$, $j=1,...,N_{F}$, that minimises
  \begin{equation}
    \sum_{i=1}^{N_{F}}\sum_{j=1}^{N_{F}}D_{i,j}\norm*{X^{i,n} - X^{j,n}}^2,
  \end{equation}
  subject to the marginal constraints
  \begin{equation}
    \sum_{i=1}^{N_{F}} D_{i,j} = \frac{1}{N_{F}},
    \quad \sum_{j=1}^{N_{F}} D_{i,j} = w^{i,n}.
  \end{equation}
  This is equivalent to maximising the covariance between them. After
 using linear programming to obtain the
 minimal matrix $D_{i,j}$, we compute the ensemble transform
  \begin{equation}
    \tilde{X}^{j,n} = N_{F} \sum_{i=1}^{N_{F}} D_{i,j}X^{i,n}, \quad j=1,\ldots,N_{F}.
  \end{equation}
  
  We verify that the mean of this new ensemble was preserved from the weighted ensemble
  
  \begin{equation}
    \frac{1}{N_{F}}\sum_{j=1}^{N_{F}}\tilde{X}^{j,n} = \sum_{j=1}^{N_{F}}\sum_{i=1}^{N_{F}} D_{i,j}X^{i,n} = \sum_{i=1}^{N_{F}} w^{i,n}X^{i,n}.
  \end{equation}
  
  It was demonstrated in \cite{Reich} that
  this transformation provides weakly converging approximations of
  higher moments as $N_{F}\to \infty$.

\smallskip
  
\item {\bfseries Re-couple the pair of ensembles in each difference
    estimator.} 
    
    For each difference estimator, given the transformed coarse and fine ensembles $\{\tilde{X}_C^{i, n}\}_{i=1}^{N_{F}}$
    and $\{\tilde{X}_F^{i, n}\}_{i=1}^{N_{F}}$ with equal weights $1/N_{F}$ respectively, we seek the
    coupling matrix $T_{i,j}$, $i=1,...,N_{F}$, $j=1,...,N_{F}$, that minimises
    \begin{equation}
      \sum_{i=1}^{N_{F}}\sum_{j=1}^{N_{F}} T_{i,j}\norm*{\tilde{X}_C^{i, n} - \tilde{X}_F^{j, n}}^2,
    \end{equation}
    where $T_{i,j}$ must take non-negative integer values subject to
    the constraints
    \begin{equation}
      \sum_{i=1}^{N_{F}} T_{i,j} = \sum_{j=1}^{N_{F}} T_{i,j} = \frac{1}{N_F}.
    \end{equation}
    This is an assignment problem that can be solved by the Hungarian
    algorithm, resulting in a matrix $T_{i,j}$ with all entries equal
    to either $1/N_F$ or $0$. We then re-order the coarse ensemble so
    that $T_{i,j}$ becomes a diagonal matrix.
    
\end{enumerate}

\subsection{The Seamless Multilevel Ensemble Transform Particle Filter}

The multilevel ensemble transform filter has an inelegant aspect, in
that we first decorrelate the pairs of ensembles through independent
ensemble transformation, before subsequently restoring correlation by
using an assignment problem. We should expect the coupling between the pairs of ensembles to be stronger if we were to avoid this aspect and keep them correlated throughout the transformation / coupling scheme.

To address this, we now describe a new multilevel filtering algorithm,
which we call the seamless multilevel ensemble transform filter (seamless MLETPF), to
correct this inelegance. This algorithm redesigns the assignment problem used to couple pairs of ensembles as two different optimal transportation problems. In this new version of the algorithm, after updating the weights of each pair of
ensembles, we perform the following steps for each pair of coarse
and fine ensembles in the difference estimators.

\bigskip

\begin{enumerate}

\item {\bfseries Find a coupling between the weighted fine and coarse ensembles.}

Generate a coupling matrix $D_{i,j}$, $i=1,...,N_{F}$, $j=1,...,N_{F}$ that minimizes
\begin{equation}
\sum^{N_{F}}_{i=1}\sum^{N_{F}}_{j=1}D_{i,j}\norm*{\hat{X}_{C}^{i, n}-X_{F}^{j, n}}^{2},
\label{equation:cost_function}
\end{equation}
subject to the marginal constraints given by
\begin{equation}
\sum^{N_{F}}_{j=1}D_{i,j}=\hat{w}_{C}^{i, n}, \qquad \sum^{N_{F}}_{i=1}D_{i,j}=w_{F}^{j, n}.
\end{equation}
We define an intermediate ensemble given by
\begin{equation}
X_{C*}^{j, n} = \sum^{N_{F}}_{i=1} D_{i,j}\hat{X}_{C}^{i, n}(w_{F}^{j, n})^{-1},
\end{equation}
for $j=1,...,N_{F}$ with importance weights $\{w_{F}^{j, n}\}_{j=1}^{N_{F}}$.

\smallskip

\item { \bfseries Transform the fine ensemble.}

Generate a coupling matrix $T_{i,j}$, $i=1,...,N_{F}$, $j=1,...,N_{F}$ that minimizes

\begin{equation}
\sum^{N_{F}}_{i=1}\sum^{N_{F}}_{j=1}T_{i,j}\norm*{X_{F}^{i, n}-X_{F}^{j, n}}^{2},
\label{equation:transform_cost_function}
\end{equation}
subject to the marginal constraints given by
\begin{equation}
\sum^{N_{F}}_{j=1}T_{i,j}=w_{F}^{i, n}, \qquad \sum^{N_{F}}_{i=1}T_{i,j}=\frac{1}{N_{F}}.
\end{equation}
This forms the matrix needed to transform the weighted ensemble $\{X_{F}^{i, n}, w_{F}^{i, n}\}_{i=1}^{N_F}$ into the evenly weighted ensemble $\{\tilde{X}_{F}^{j, n}\}_{j=1}^{N_F}$, where $\tilde{X}_{F}^{j, n}=\sum^{N_{F}}_{i=1}T_{i,j}N_{F}X_{F}^{i, n}$.

\smallskip

\item { \bfseries Transform the coarse ensemble with fine weights to an evenly weighted ensemble.}

Find another coupling matrix $\tilde{T}_{i,j}$, $i=1,...,N_{F}$, $j=1,...,N_{F}$ that minimizes

\begin{equation}
\sum^{N_{F}}_{i=1}\sum^{N_{F}}_{j=1}\tilde{T}_{i,j}\norm*{X_{C*}^{i, n}-\tilde{X}_{F}^{j, n}}^{2},
\end{equation}
subject to the marginal constaints given by
\begin{equation}
\sum^{N_{F}}_{j=1}\tilde{T}_{i,j}=w_{F}^{i, n}, \qquad \sum^{N_{F}}_{i=1}\tilde{T}_{i,j}=\frac{1}{N_{F}}.
\end{equation}
Here we use the evenly weighted new transformed finer ensemble $\{\tilde{X}_{F}^{i, n}\}_{i=1}^{N_{F}}$ in the cost function for $\tilde{T}_{i,j}$ to keep the coarse and fine ensembles closely coupled. Finally the new transformed, evenly weighted coarse ensemble $\{\tilde{X}_{C}^{j, n}\}_{j=1}^{N_{F}}$ is given by
\begin{equation}
\tilde{X}_{C}^{j, n}=\sum^{N_{F}}_{i=1}\sum^{N_{F}}_{k=1}\tilde{T}_{i,j}N_{F}D_{k,i}\big(w_{F}^{i, n}\big)^{-1}\hat{X}_{C}^{k, n},
\end{equation}
for $j=1,...,N_{F}$.

\end{enumerate}

\bigskip

As mentioned previously, this method aims to couple the ensembles in each multilevel difference estimator, so that the sample covariance of the difference between them, $\mathbb{V}_{l}^n=\Cov[\tilde{X}_{F}^n-\tilde{X}_{C}^n]$, decays asymptotically as $l \to \infty$ at a rate $\mathcal{O}(h_{l}^{\beta})$. Given other assumptions, \cite{Cliffe} shows that if $\beta>\gamma$, where $\gamma$ is defined as the model cost exponent earlier in this paper, an optimal computational cost reduction can be reached. The value of $\beta$ that this seamless coupling achieves appears, from numerical studies at the end of the paper, to be greater than the previous methodology outlined in \cite{Gregory} for multivariate problems. Greater values of $\beta$, seemingly offered by this new algorithm, therefore signify the potential of greater computational cost reductions in some cases.

The assignment problem in the algorithm in \cite{Gregory} could be interpreted as rearranging the transformed particles in both the coarse and fine ensembles in each multilevel difference estimator. By replacing this problem with the seamless combination of optimal transportation problems, we instead modify the particles to maintain as strong a coupling as possible between the two ensembles, throughout each step of the algorithm. Thus we are optimising over a continuous space of couplings instead of a discrete space of couplings.

One can verify that the sample means of the coupled, transformed $\{\tilde{X}_{C}^{i, n}\}_{i=1}^{N_{F}}$ and $\{\tilde{X}_{F}^{i, n}\}_{i=1}^{N_{F}}$ preserve those from the weighted ensembles given before the transform by
\begin{equation}
\frac{1}{N_{F}}\sum^{N_{F}}_{j=1}\tilde{X}_{F}^{j, n}=\sum^{N_{F}}_{j=1}\sum^{N_{F}}_{i=1}\frac{1}{N_{F}}T_{i,j}N_{F}X_{F}^{i, n}=\sum^{N_{F}}_{i=1}w_{F}^{i, n}X_{F}^{i, n},
\end{equation}
for the finer level \cite{Reich}, and
\begin{equation}
\begin{split}
\frac{1}{N_{F}}\sum^{N_{F}}_{j=1}\tilde{X}_{C}^{j, n}&=\sum^{N_{F}}_{j=1}\sum^{N_{F}}_{i=1}\sum^{N_{F}}_{k=1}\frac{1}{N_{F}}\tilde{T}_{i,j}N_{F}D_{k,i}\big(w_{F}^{i, n}\big)^{-1}\hat{X}_{C}^{k, n}\\
\quad &=\sum^{N_{F}}_{i=1}\sum^{N_{F}}_{k=1}D_{k,i}\hat{X}_{C}^{k, n}=\sum^{N_{F}}_{k=1}\hat{w}_{C}^{k, n}\hat{X}_{C}^{k, n},
\end{split}
\end{equation}
for the coarse level. In terms of studying the asymptotic consistency of the higher moments of $\{\tilde{X}_{C}^{i, n}\}_{i=1}^{N_{F}}$, we consider the coarse transform as a combination of two ensemble transforms. \cite{Reich} showed that the first transform generating the intermediate ensemble will satisfy weak asymptotic convergence
\begin{equation}
\sum^{N_F}_{i=1}g(X_{C*}^{i,n})w_{F}^{i,n} \to \sum^{N_F}_{i=1}g(\hat{X}_{C}^{i,n})\hat{w}_{C}^{i,n} ,
\end{equation}
as $N_F \to \infty$, and by the same logic
\begin{equation}
\frac{1}{N_F}\sum^{N_F}_{i=1}g(\tilde{X}_{C}^{i,n}) \to \sum^{N_F}_{i=1}g(X_{C*}^{i,n})w_{F}^{i,n} ,
\end{equation}
also as $N_F \to \infty$. Thus
\begin{equation}
\frac{1}{N_F}\sum^{N_F}_{i=1}g(\tilde{X}_{C}^{i,n}) \to \sum^{N_F}_{i=1}g(\hat{X}_{C}^{i,n})\hat{w}_{C}^{i,n} ,
\end{equation}
as $N_F \to \infty$. We demonstrate this numerically in the next section. Thus the coupled, transformed coarse ensemble in each of the difference estimators in the seamless MLETPF produce consistent estimators to statistics of $X_{C}^n$ and therefore the telescoping formula in Section \ref{subsec:ml_filter} is asymptotically consistent.

One can also use localisation to reduce the dimensionality of the cost functions (and to reduce the computational cost of the linear transport problems) in the same way as was mentioned earlier in the paper, by simply implementing the above algorithm for every individual component of a multivariate random variable. 

When localisation is not employed, the forward model cost of the standard ETPF (with a fixed order of accuracy) must dominate over the optimal transportation cost of the seamless MLETPF in order for the overall computational cost of the ETPF to be reduced. Suppose we desire a mean-square-error (MSE) of $\mathcal{O}(\epsilon^{2})$, for a small $\epsilon$; we require $N_{0}=\mathcal{O}(\epsilon^{-2})$ (for the seamless MLETPF) and $N=\mathcal{O}(\epsilon^{-2})$ (for the standard ETPF). Also assume $h_{L} = \mathcal{O}(\epsilon)$ for a first order accurate discretization technique. Then the forward model cost of the standard ETPF is
$C_{FM} = c_{m}\left(\epsilon^{-2-\gamma}\right)$,
where $c_{m}$ is a constant.

If we assume that $2N_{l}=N_{l-1}$ for simplicity, the optimal transportation cost of the seamless MLETPF can be written as
\begin{equation}
\begin{split}
C_{OT}&=c_{t}N_{0}^{3}\sum^{L}_{l=0}\left(1/2\right)^{3l}\log\left(N_{0}\left(1/2\right)^{l}\right) \leq (8/7)c_{t}\epsilon^{-6}\log\left(\epsilon^{-2}(1/2)^{1/7}\right),
\end{split}
\end{equation}
using standard (arithmetico) geometric series results. Here $c_{t}$ is a constant. In the case where $\gamma$ is not extremely large, $C_{FM}$ will become less than $C_{OT}$ as $\epsilon \to 0$, and the optimal transportation cost will dominate; forward model cost speed-ups from the multilevel framework would be worthless. However we expect that the seamless MLETPF will offer overall speed-ups in a certain $\epsilon$-regime, provided that $c_{m}/c_{t} \gg 1$.

\section{Numerical Examples}

The three numerical experiments used in this section are framed to test out the posterior consistency of the seamless MLETPF, the improved rate of variance decay for a non-localised filter estimate from the algorithm in \cite{Gregory}, and the computational cost reductions for a fixed accuracy from localised ETPF estimators.
In some of the numerical examples we need to estimate the root-mean-square-error (RMSE) of the MLETPF; we take this relative to an accurate (with variance / discretization bias much lower than any of the experimental MLETPF simulations) ETPF approximation. The details of the discretization level and sample size of this approximation will be given during these examples. The RMSE, against $E[X^{n}]$ conditioned on the observations, or an accurate approximation of it as the case here may be, is estimated over the time indices $n\in[1,N_{y}]$ via
\begin{equation}
RMSE = \sqrt{\frac{1}{N_{y}}\sum^{N_{y}}_{n=1}\norm*{\hat{E}^{N_{0},...,N_{L}}[X_{L}^{n}] - E[X^{n}]}^2}.
\end{equation}
Computational cost is here given as a runtime in seconds using a standard Python implementation of the algorithms.

\bigskip

\begin{figure}[ht]
%\hspace{\fill}
%\begin{minipage}{.45\textwidth}
  \centering
  \includegraphics[width=125mm]{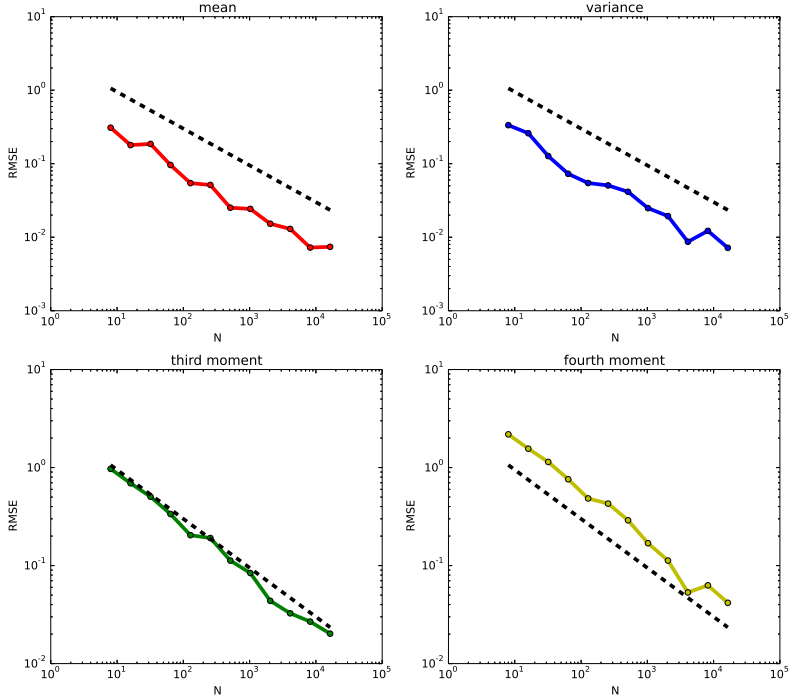}
  \caption{\textit{Root-mean-square-error (from the mean, variance and estimated third and fourth moments of the true posterior) of the sample mean, variance, third and fourth moments of $\tilde{X}_{C}$, the posterior of $X_{C}$ given an observation, generated from the seamless transform / coupling algorithm on the ensembles, size $N$, sampled from $X_{C}$, $X_{F}$. The black dashed lines represent square root convergence.}}
  \label{figure:consistency_higher_moments_seamless_coupling}
%\end{minipage}
\end{figure}

\noindent \textit{Example (Consistency of posterior statistics).} Here, the consistency of the mean, variance, third and fourth moments of the coarse posterior after a single seamless transform / coupling step will be confirmed numerically. As the fine ensemble is transformed using the standard ensemble transform method, one notes these moments of this posterior are consistent via the original literature in \cite{Reich}. We use
\begin{equation}
X_{C} \sim N\big(1,1\big), \qquad X_{F} \sim N\big(0.5,1\big).
\end{equation}
The single observation is given by $Y_{obs}=0.1$ with a measurement error variance of $2$. Using Bayes' Theorem, the true coarse posterior, given the observation, and the prior $X_{C}$, is $N\big(0.7,2/3\big)$. Ensembles of varying size $N$ are drawn from $X_{C}$ and $X_{F}$, weighted with respect to the observation distribution and finally $\tilde{X}_{C}$ is found using the seamless coupling / transform algorithm. Figure \ref{figure:consistency_higher_moments_seamless_coupling} shows the root-mean-square-error of the sample mean, variance, third and fourth moments of $\tilde{X}_{C}$ from that of the true posterior over 10 independent ensembles.

\bigskip

\noindent {\textit{Example (Lorenz 63 equations).} A stochastic Lorenz-63 system is now used to test the seamless ensemble transform / coupling algorithm on a multivariate problem. The equations are given by the 3-component chaotic nonlinear system in $X=(x,y,z)$
\begin{equation}
\frac{dX}{dt}=
\begin{cases}
\sigma(y-x)+\phi\frac{dW}{dt},& \\
x(\rho-z)-y+\phi\frac{dW}{dt},& \\
xy-\beta z+\phi\frac{dW}{dt},&
\end{cases}
\end{equation}
with $\rho=28$, $\sigma=10$, $\beta=8/3$, $\phi=0.1$. The scalar Brownian motion $W$ in the above system will be the same for each component to maximise the impact of the strong nonlinearity in the equations, i.e. to make it a more challenging filtering example.
Let $X^n$ be the solution to the above Lorenz-63 equations at time $t_{n}$, and let $X_{l}^{n}$ be the discretization of $X^{n}$ using the forward Euler numerical scheme with time-step $h_{l}=2^{-9-l}$.
The seamless MLETPF estimator with $N_{l}=2^{8-l}$, $l \in [0,6]$, using the full optimal transport problems with no localisation, is used to compute an approximation to $E[X^{n}]$. Here $\Delta t=2^2h_{0}$ and we take $n \in [1,1280]$. This coarsest time-step of $h_{0}=2^{-9}$ is used due to stability reasons. Observations are given by a measurement error with $R=0.25I$, where $I$ is the $3 \times 3$ identity matrix.

First, the example will be used to show to the improvement in the rate of variance decay of the coupling between fine and coarse ensembles with the level of hierarchy, $\mathbb{V}_{l}^n=\mathcal{O}(h_{l}^{\beta})$, in a multivariate case. Figure \ref{figure:Lorenz63SeamlessCoupling} show the average rates of asymptotic decay of the sample covariance, from using 5 independent implementations of the seamless MLETPF estimator and the original algorithm in \cite{Gregory}. We observe that the seamless algorithm now produces a variance decay close to $\beta = 2$. The rate is certainly improved from the original methodology, in \cite{Gregory}, where $\beta\approx 1$ is obtained. 

Second, this example shows the effects that the increased rate of variance decay, shown in Figure \ref{figure:Lorenz63SeamlessCoupling}, has on the overall forward model cost of the seamless MLETPF estimator. Given that optimal transportation algorithmic costs dominate over the cheap forward model ($\gamma=1$) in this example, there isn't an overall speed-up in convergence offered by either implementation of the MLETPF from that of the standard ETPF. However the increased rate of variance decay in the seamless MLETPF estimator can be seen to produce further reductions in the forward model costs from that of the standard MLETPF and ETPF estimators. This means that there is an overall runtime benefit whenever the forward model cost of the standard ETPF dominates the optimal transportation costs.

We define a desired order of magnitude of RMSE, $\epsilon$, and use this to determine the parameters $L$ and $N_{l}$ as done in \cite{Gregory}. From this, we set $L=\ceil*{-log(\epsilon^{-1})/log(2)}$ and $N_{l}=\ceil*{\epsilon^{-2}2^{-l}}$ for the standard implementation of the MLETPF and $N_{l}=\ceil*{\epsilon^{-2}2^{-(3/2)l}}$ for the seamless one. The choice of using the geometric decay in $N_{l}$ is primarily for simplicity and a proof of concept; we do not claim these are optimal as would be obtained by using the algorithm in \cite{Giles}. The different rates of geometric decay in $N_{l}$ are designed to exploit the different rates of sample variance decay shown in Figure \ref{figure:Lorenz63SeamlessCoupling}, in order to gain a reduction in the growth of forward model cost for the seamless MLETPF for a fixed order of estimator variance $\left(\sum^{L}_{l=0}\mathbb{V}_{l}/N_{l}\right)$, and thus accuracy. For the standard ETPF, $N=\ceil*{\epsilon^{-2}}$ and $L$ is set to be the same as in the MLETPF estimators. The reduction in the growth of forward model cost for the seamless MLETPF is shown in Figure \ref{figure:Lorenz63SeamlessCouplingCompCost}, for decreasing values of $\epsilon$. The rates, $\mathcal{O}\left(\epsilon^{-2}\right)$ for the seamless implementation and $\mathcal{O}\left(\epsilon^{-2}\log(\epsilon)^{2}\right)$ for the standard implementation, are consistent with the analysis in \cite{Giles} for when $\beta>\gamma$ and $\beta=\gamma$ respectively. Here computational cost is measured by the total number of floating point operations for the forward model computation.

\begin{figure}[ht]
%\hspace{\fill}
%\begin{minipage}{.45\textwidth}
  \centering
  \includegraphics[width=120mm]{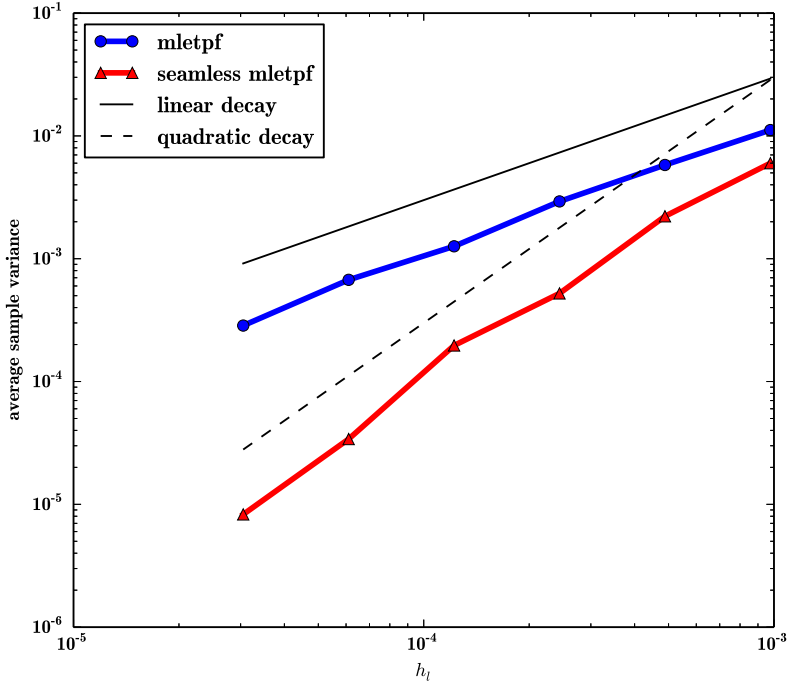}
  \caption{\textit{Average, over 5 independent simulations and all assimilation steps $n \in [1,1280]$, of the sample variance $Tr\big(\mathbb{V}_{l}^n\big)=Tr\big(\Cov[\tilde{X}_{l}^n-\tilde{X}_{l}^n]\big)$ for $l \in [1,6]$. Here $X$ denotes the solution to the Lorenz 63 system. Results are shown for the seamless MLETPF and the previous algorithm of \cite{Gregory} in the red (with triangular points) and blue (with circular points) lines respectively. Asymptotes of linear and quadratic decay, with decreasing $h_{l}$, are shown by black solid and dashed lines respectively.}}
  \label{figure:Lorenz63SeamlessCoupling}
%\end{minipage}
\end{figure}

\begin{figure}[ht]
%\hspace{\fill}
%\begin{minipage}{.45\textwidth}
  \centering
  \includegraphics[width=135mm]{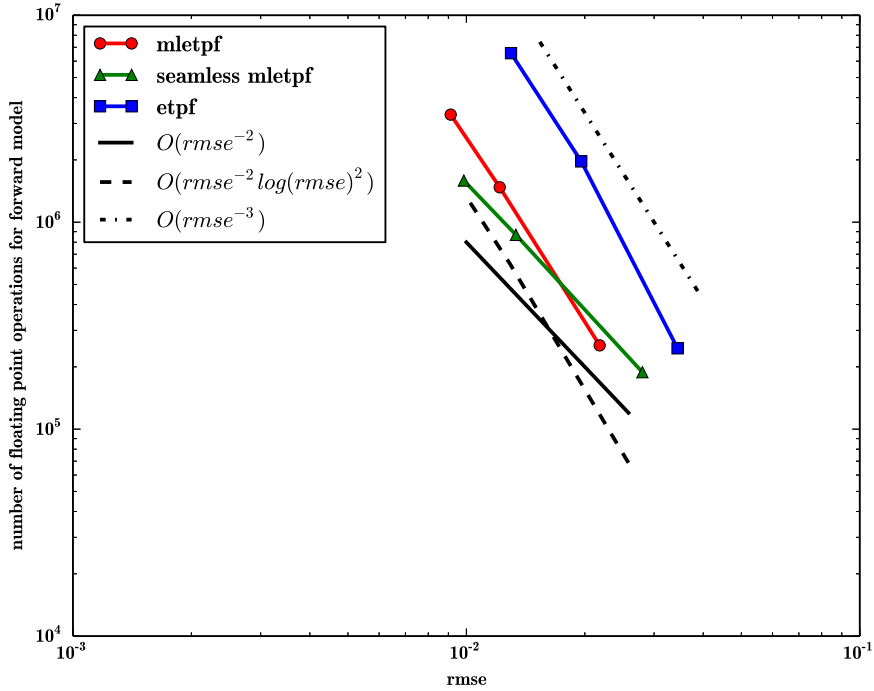}
  \caption{\textit{Forward model computational costs against RMSE for decreasing $\epsilon$ values,  of the seamless and standard MLETPF estimators and the standard ETPF (without localisation) in the Lorenz 63 problem set-up.}}
  \label{figure:Lorenz63SeamlessCouplingCompCost}
%\end{minipage}
\end{figure}

\bigskip

\noindent \textit{Example (Lorenz 96 equations).} The Lorenz-96 system is given by the following spatial discretization,

\begin{equation}
\frac{dX(j)}{dt}=-\frac{\big(X(j-1)X(j+1)-X(j-2)X(j-1)\big)}{3\Delta}-X(j)+F+\sigma^2\frac{dW_j}{dt}.
\end{equation}

Here, $F$ is a forcing term, typically taken to be $F=8$ for chaotic behaviour, $j \in [1,40]$, $\Delta = 0.5$, and $W_{j}$ are independent standard Brownian Motions, with $\sigma^{2} = 0.1$. This independence allows the reduction in spatial correlations in this system and hence localisation, which will be used for this problem, has little effect on the performance of the assimilation, as in \cite{Gregory}.
For systems with strong spatial dependence structures, localisation will have more of an affect on the ETPF and it's multilevel counterparts \cite{Cheng}. The choice of some parameters, such as the frequency of assimilation, is an important aspect to consider when choosing if localisation, and thus the ETPF is suitable for a particular problem. For example, if assimilation is more frequent, spatial dependence will be altered more frequently, that in turn could lead to a worse performing filter.

Let $X_{l}^n(j)$ represent the $j$'th component of $X$ at $t_{n}$ (with time-index $n$), discretized by the forward Euler numerical scheme with time-step $h_{l}=2^{-8-l}$. The coarsest time-step, $h_{0}=2^{-8}$, is used here for stability purposes. Here, $\Delta t = h_{0}$, and we take $n \in [1, 1280]$; observations are given by using a measurement error variance of $R=0.25I$, where $I$ is the $40 \times 40$ identity matrix.
This system will be used to test out the computational cost reductions of the seamless MLETPF, from that of the standard ETPF, and the standard MLETPF implementation, in finding estimates to statistics of $X^n$ of a pre-defined accuracy. Due to the use of localisation, the model cost will now dominate over the cost of the ensemble transform  / coupling scheme and so we evaluate whether the seamless MLETPF offers overall computational cost reductions. As in \cite{Gregory}, the use of localisation does lead to an inconsistent approximation to statistics of $X^n$ so we measure RMSE relative to a localised, high accuracy ETPF approximation, that the estimates are asymptotically consistent with, as stated at the start of this section. The discretization level of this approximation is $L=13$ and the sample size is $N=40000$.

The authors assume that one has already chosen to use the fully localised ETPF for this problem and are using the multilevel framework to improve the efficiency, so hence evaluating the convergence of the MLETPF to a localised ETPF approximation rather than the truth is sensible in this setting. If one were to compare to the truth here, the speed-up in convergence that would be seen in the multilevel cases would plateau due to the localisation bias; this is a problem inherent of localisation in general and not specifically of the multilevel framework.

Pre-defined values of $L$ and $N_{l}$ will be used in the seamless and standard MLETPF estimators for a user-defined $\epsilon$ which determines an order of magnitude of RMSE as done in the previous example. We set $L=\ceil*{-\log(\epsilon)/\log(2)}$ and $N_{l}=\ceil*{\epsilon^{-2}2^{-(3/2)l}}$ for each $\epsilon$ respectively in both the seamless and standard implementations of the MLETPF. For the standard ETPF estimator, we set the sample size to be $N=\ceil*{\epsilon^{-2}}$ and $L$ to be the same as above.
In Figure \ref{figure:seamless_coupling_variance_lorenz96}, the average decay of sample variance, $\mathbb{V}_{l}$, with increasing $l$ is shown for both the standard and seamless implementations of the MLETPF. As was observed in \cite{Gregory}, the standard implementation of the MLETPF achieves the same optimal quadratic rate of variance decay (as the seamless one) in this fully localised case, however the magnitude of variance is smaller for the seamless implementation.

The RMSE against computational cost for the standard and seamless implementations of the MLETPF estimator, as well as the standard ETPF estimator, approximating $E[X^n]$, are shown for different $\epsilon$ values in Figure \ref{figure:seamless_coupling_cost_accuracy_lorenz96}. One notes that in this problem, where $\gamma=1$, the computational cost of the standard ETPF estimator follows $\mathcal{O}(\epsilon^{-3})$ scaling. Computational cost reductions, down to approximately $\mathcal{O}(\epsilon^{-2})$ scaling, are seen for both the standard and seamless implementations of the MLETPF estimator. However the computational costs for the seamless MLETPF estimator are of lower magnitude along this scaling than in the standard implementation of the MLETPF. Similar rates of convergence are shown for the three estimators approximating $E\left[\left(X^n\right)^{2}\right]$ in Figure \ref{figure:seamless_coupling_cost_accuracy_lorenz96_second_moment}.

\begin{figure}[ht]
%\hspace{\fill}
%\begin{minipage}{.45\textwidth}
  \centering
  \includegraphics[width=135mm]{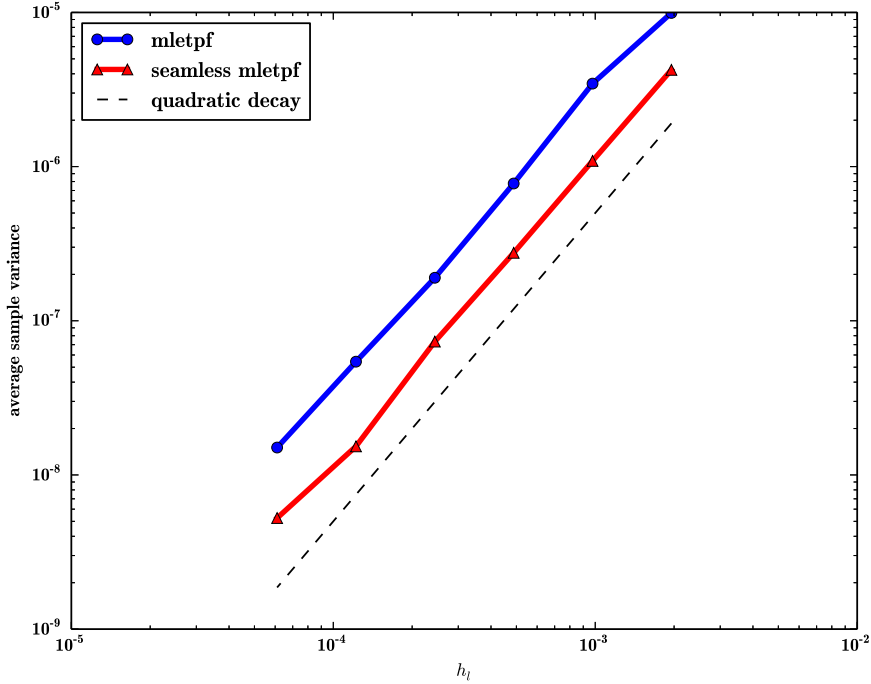}
  \caption{\textit{Average, over 5 independent simulations and all assimilation steps $n \in [1,1280]$, of the sample variance $Tr\big(\mathbb{V}_{l}^n\big)=Tr\big(\Cov[\tilde{X}_{l}^n-\tilde{X}_{l}^n]\big)$ for $l \in [1,6]$. Here $X$ denotes the solution to the Lorenz 96 system. Results are shown for the seamless MLETPF and the previous algorithm of \cite{Gregory} (both localised) in the red (with triangular points) and blue (with circular points) lines respectively. An asymptote of quadratic decay, with decreasing $h_{l}$, is shown by the dashed line.}}
   \label{figure:seamless_coupling_variance_lorenz96}
%\end{minipage}
\end{figure}

\begin{figure}[ht]
%\hspace{\fill}
%\begin{minipage}{.45\textwidth}
  \centering
  \includegraphics[width=135mm]{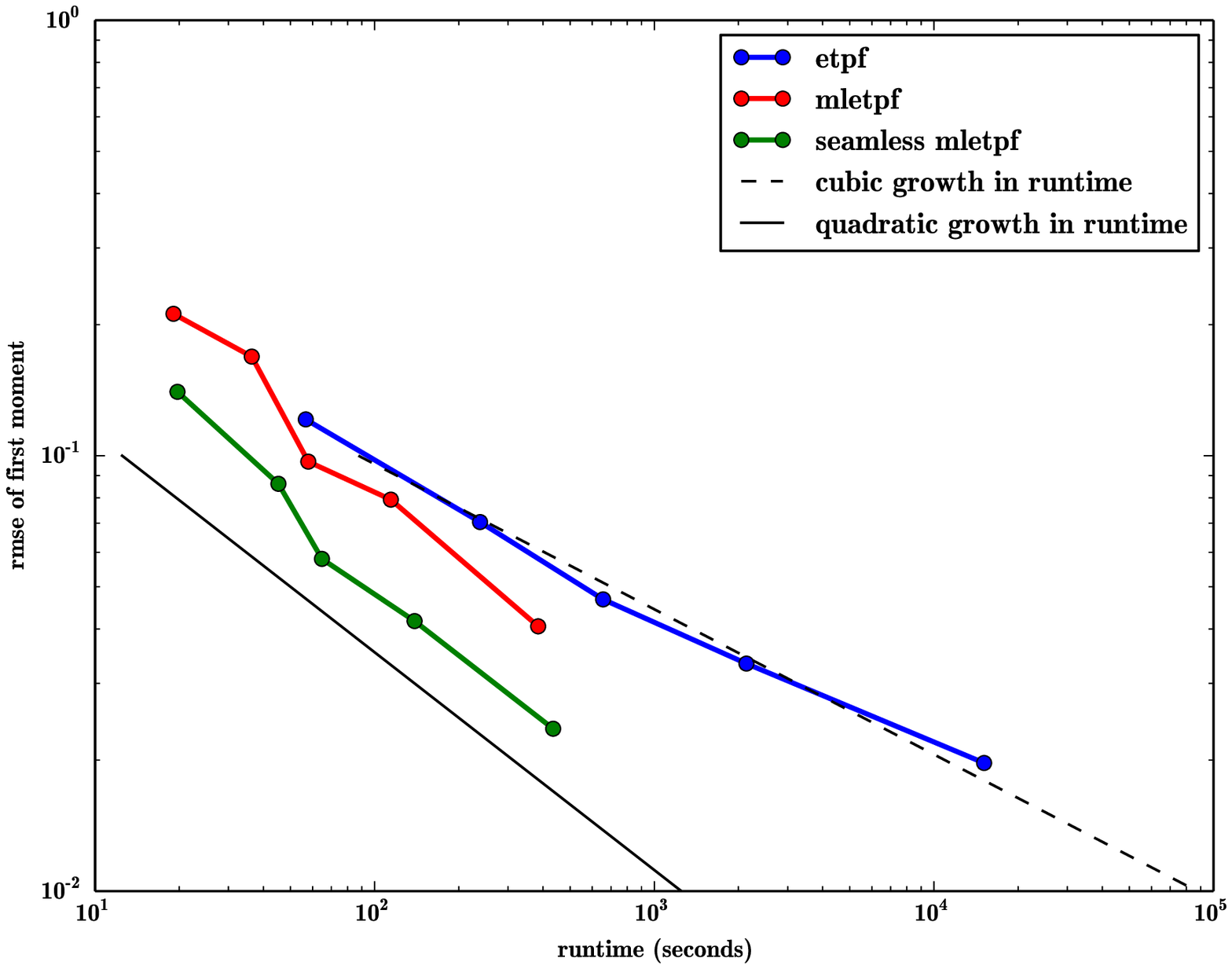}
  \caption{\textit{Root-mean-square-error (RMSE) against runtime of the approximations to $E\left[X^{n}\right]$ from the standard ETPF and both implementations of the MLETPF for decreasing values of $\epsilon$ in the Lorenz 96 problem set-up. Rates of $\mathcal{O}(\epsilon^{-3})$ and $\mathcal{O}(\epsilon^{-2})$ of computational cost growth are shown in black dashed and solid lines respectively. }}
  \label{figure:seamless_coupling_cost_accuracy_lorenz96}
%\end{minipage}
\end{figure}

\begin{figure}[ht]
%\hspace{\fill}
%\begin{minipage}{.45\textwidth}
  \centering
  \includegraphics[width=135mm]{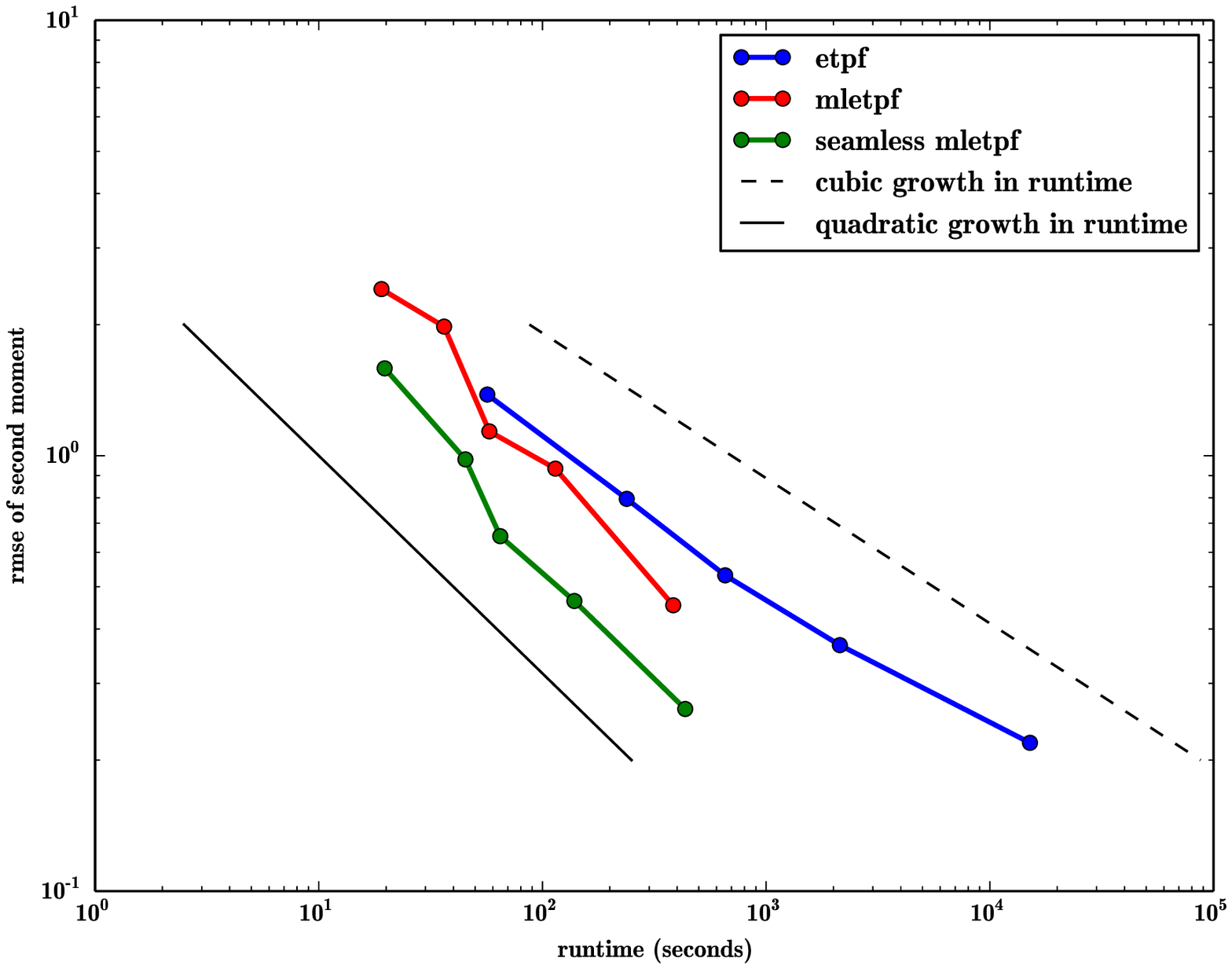}
  \caption{\textit{Root-mean-square-error (RMSE) against runtime of the approximations to $E\left[\left(X^{n}\right)^{2}\right]$ from the standard ETPF and both implementations of the MLETPF for decreasing values of $\epsilon$ in the Lorenz 96 problem set-up. Rates of $\mathcal{O}(\epsilon^{-3})$ and $\mathcal{O}(\epsilon^{-2})$ of computational cost growth are shown in black dashed and solid lines respectively. }}
  \label{figure:seamless_coupling_cost_accuracy_lorenz96_second_moment}
%\end{minipage}
\end{figure}

\section{Discussion and Summary}

This paper has presented a seamless algorithm for an efficient version of the Ensemble Transform Particle Filter (ETPF), the multilevel ETPF (MLETPF), improving on an algorithm in \cite{Gregory}. It is seamless in the sense that we use a combination of three optimal transport problems to transform and positively couple coarse and fine ensembles simultaneously. This replaced the assignment problem that was solved to recouple them, after transforming them each independently, in the aforementioned previous algorithm.

The benefit of this change is that, from a proof of concept presented in this paper, it can be seen that using the seamless MLETPF, reduces the variance between the coarse and fine ensembles, and even the rate at which this variance decays with increasing levels of accuracy in some cases, from that of \cite{Gregory}. This in turn leads to a lower overall variance of the MLETPF estimators, and thus can reduce the forward model cost of computing them for a fixed accuracy.

For cases where the forward model cost is low, the optimal transport cost is the computational bottleneck in both the ETPF and the multilevel equivalent unless localisation is used. Iterative methods, such as the Sinkhorn approximation in \cite{Sen}, can be used to reduce this algorithmic cost. They have been utilised in \cite{DeWiljes} by using a second-order accurate framework for the ETPF. This is an important future direction of this multilevel research, since the proposed algorithm in this paper involves a combination of optimal transport problems that could be solved using such methods.

\bibliographystyle{siam}
\bibliography{refs}

\end{document}